\documentclass[12]{article}
\usepackage{fleqn}
\usepackage{graphicx}
\begin{document}
\Large
\begin{center}{\bf
Projective Planes Over ``Galois" Double
Numbers and a Geometrical Principle of Complementarity}
\end{center}
\vspace*{.4cm}
\begin{center}
Metod Saniga$^{\dag}$ and Michel Planat$^{\ddag}$
\end{center}
\vspace*{.2cm} \normalsize
\begin{center}
$^{\dag}$Astronomical Institute, Slovak Academy of Sciences,
SK-05960 Tatransk\' a Lomnica\\ Slovak Republic\\
(msaniga@astro.sk)

\vspace*{.1cm}
 and

\vspace*{.1cm} $^{\ddag}$Institut FEMTO-ST, CNRS, D\' epartement LPMO, 32 Avenue de
l'Observatoire\\ F-25044 Besan\c con, France\\
(planat@lpmo.edu)
\end{center}

\vspace*{.0cm} \noindent \hrulefill

\vspace*{.1cm} \noindent {\bf Abstract}

\noindent
The paper deals with a particular type of a projective ring plane defined
over the ring of double numbers over Galois fields,
$R_{\otimes}(q) \equiv $ GF($q$) $\otimes$ GF($q$) $\cong$ GF($q$)[$x$]/($x(x-1$)). The plane is endowed with $(q^2 + q + 1)^2$ points/lines
and there are $(q + 1)^2$ points/lines incident with any line/point. As $R_{\otimes}(q)$ features two maximal ideals,
the neighbour relation is
not an equivalence relation, i.\,e. the sets of neighbour points to two distant points overlap. Given a point of the plane, there are $2q(q+1)$
neighbour points to it. These form two disjoint, equally-populated families under the reduction modulo either of the ideals. The points of the first
family merge with (the image of) the point in question, while the points of the other family go in a one-to-one fashion to the remaining $q(q + 1)$ points of
the associated ordinary (Galois) projective plane of order $q$. The families swap their roles when switching from one ideal to the other, which
can be regarded as a remarkable, finite algebraic geometrical manifestation/representation of the principle of complementarity. Possible domains of
application of this finding in (quantum) physics, physical chemistry and neurophysiology are briefly mentioned.\\



\noindent {\bf Keywords:} Projective Ring Planes --- Rings of
Double Numbers Over Galois Fields --- \\
\hspace*{2.05cm}Neighbour/Distant Relation --- Geometrical
Complementarity Principle

\vspace*{-.1cm} \noindent \hrulefill

\vspace*{.3cm}  \noindent
\section{Introduction}
Although (finite) projective ring planes represent a well-studied, important and venerable branch of algebraic geometry \cite{vk81}--\cite{vk95} and are
endowed with a number of fascinating, although rather counter-intuitive properties not having analogues in ordinary (Galois) projective
planes, it may well  come as a big surprise that, as far as we know, they have so far successfully evaded the attention of physicists
and/or scholars of other natural sciences. The only exception in this respect seems to be our very recent paper \cite{spl} in which we pointed out the
importance of the structure of perhaps the best known of finite ring planes, that of Hjelmslev \cite{h}--\cite{dj}, for getting deeper insights into
the properties of finite dimensional Hilbert spaces of quantum (information) theory. This paper aims at examining another interesting type
of finite projective planes, namely that defined over the ring of double numbers over a Galois field. As this ring is not a local ring like
those that serve as coordinates for Hjelmslev planes, the structure of the corresponding plane is more intricate when compared
with the corresponding Hjelmslev case and, as we shall see, may thus lend itself to more intriguing potential applications.

\section{Rudiments of Ring Theory}
In this section we recollect some basic definitions and properties
of rings that will be employed in the sequel and to the extent
that even the reader not well-versed in the ring theory should be
able to follow the paper without the urgent need of consulting
further relevant literature (e.g., \cite{fr}--\cite{ra}).

A {\it ring} is a set $R$ (or, more specifically, ($R,+,*$)) with
two binary operations, usually called addition ($+$) and
multiplication ($*$), such that $R$ is an abelian group under
addition and a semigroup under
multiplication, with multiplication being both left and right
distributive over addition.\footnote{It is customary to denote
multiplication in  a ring simply by juxtaposition, using $ab$ in
place of $a*b$, and we shall follow this convention.} A ring in
which the multiplication is commutative is a commutative ring. A
ring $R$ with a multiplicative identity 1 such that 1$r$ = $r$1 = $r$
for all $r \in R$ is a ring with unity. A ring containing a finite
number of elements is a finite ring. In what follows the word ring
will always mean a commutative ring with unity.

An element $r$ of the ring $R$ is a {\it unit} (or an invertible element) if there exists an element $r^{-1}$ such that $rr^{-1} = r^{-1} r=1$.
This element, uniquely determined by $r$,  is called the multiplicative inverse of $r$. The set of units forms a group under multiplication.
An element $r$ of $R$ is said to be a {\it zero-divisor} if there exists $s \neq 0$ such that $sr= rs=0$. In a finite ring, an element is either a
unit or a  zero-divisor. A ring in which every non-zero element is a unit is a {\it field}; finite (or Galois) fields, often denoted by GF($q$),
have $q$ elements and exist only for $q = p^{n}$, where $p$ is a prime number and $n$ a positive integer.

An {\it ideal} ${\cal I}$ of $R$ is a subgroup of $(R,+)$ such
that $a{\cal I} = {\cal I}a \subseteq {\cal I}$ for all $a \in R$. An ideal of the ring $R$ which is not contained in any other ideal but $R$ itself is called a
{\it maximal}  ideal. If an ideal is of the form $Ra$ for some element $a$ of $R$ it is called a {\it principal} ideal, usually denoted by $\langle a \rangle$.
A ring with a unique maximal ideal is a {\it local} ring. Let $R$ be a ring and ${\cal I}$ one of its ideals. Then $\overline{R} \equiv R/{\cal I}
= \{a + {\cal I} ~|~ a \in R\}$ together with
addition $(a + {\cal I}) + (b + {\cal I}) = a + b +  {\cal I}$ and multiplication $(a + {\cal I})(b + {\cal I}) = ab +  {\cal I}$ is a ring, called the
quotient, or factor, ring of $R$ with respect to ${\cal I}$; if ${\cal I}$ is maximal,
then $\overline{R}$ is a field.

A mapping $\pi$:~ $R \mapsto S$ between two rings $(R,+,*)$ and
$(S,\oplus, \otimes)$ is a ring {\it homomorphism} if it meets the
following constraints: $\pi (a + b) = \pi (a) \oplus \pi (b)$,
$\pi (a * b) = \pi (a) \otimes \pi(b)$ and $\pi (1) = 1$  for any
two elements $a$ and $b$ of $R$. From this definition it is
readily discerned  that $\pi(0) = 0$, $\pi(-a) = -\pi(a)$, a unit
of $R$ is sent into a unit of $S$ and the set of elements $\{a \in
R~ |~ \pi(a) = 0\}$, called the {\it kernel} of $\pi$, is an ideal
of $R$. A bijective ring homomorphism is called a ring {\it
iso}morphism; two rings $R$ and $S$ are called isomorphic, denoted
by $R \cong S$, if there exists a ring isomorphism between them.

Finally, we mention a couple of relevant examples of rings: a polynomial ring, $R[x]$, viz. the set of all polynomials in one
variable $x$ and with coefficients in the ring $R$, and the ring $R_{\otimes}$ that is a (finite) direct product of rings,
$R_{\otimes} \equiv R_{1} \otimes R_{2} \otimes \ldots \otimes R_{n}$, where both addition and multiplication are carried out componentwise and
the component rings need not be the same.

\section{Projective Plane over the Ring of Double Numbers over a Galois Field}

The principal objective of this section is to introduce the basic properties of the projective plane defined over the direct product of two
Galois fields, $R_{\otimes}(q) \equiv$ GF($q$) $\otimes$ GF($q$) = $\{$[$a, b$]; $a, b\in$ GF($q$)$\}$ with componentwise addition
and multiplication, which is the ring isomorphic to the following quotient ring:
\begin{equation}
    R_{\otimes}(q) \cong {\rm GF}(q)[x]/(x^{2} -  x) \cong {\rm GF}(q) \oplus  e {\rm GF}(q),~e^{2} = e,~e \neq 0, 1.
\end{equation}
From the last equation it is straightforward to see that $R_{\otimes}(q)$ contains $\#_{{\rm t}} = q^{2}$ elements,
out of which there are $\#_{{\rm z}} = 2q - 1$ zero-divisors and, as the ring is finite, $\#_{{\rm u}} = \#_{{\rm t}} - \#_{{\rm z}}
= q^{2} - 2q + 1 =  (q - 1)^{2}$ units.
The set of zero-divisors consists of two maximal (and principal as well) ideals,
\begin{equation}
    {\cal I}_{\langle e \rangle} \equiv \langle e \rangle,
\end{equation}
and
\begin{equation}
    {\cal I}_{ \langle e-1 \rangle} \equiv \langle e-1 \rangle,
\end{equation}
each of cardinality $q$; we note that ``0" is the (only) common element of them.

The projective plane over $R_{\otimes}(q)$, henceforth denoted as PR$_{\otimes}(2, q)$, is a particular representative of a rich and variegated
class of projective planes defined over rings of stable rank two \cite{vk81}--\cite{vk84}, \cite{vk95}.
It is an incidence structure whose points are classes of
ordered triples
$(\varrho \breve{x}_{1}, \varrho  \breve{x}_{2}, \varrho \breve{x}_{3})$ where $\varrho$ is a unit and not all the three $\breve{x}_{i}$, if zero-divisors,
belong to the same ideal
and whose lines are, dually, ordered triples $(\varsigma \breve{l}_{1}, \varsigma \breve{l}_{2}, \varsigma \breve{l}_{3})$ with $\varsigma$ and
$\breve{l}_{i}$
enjoying the
same properties as $\varrho$ and $\breve{x}_{i}$, respectively, and where the incidence relation is defined by
\begin{equation}
\sum_{i=1}^{3}  \breve{l}_{i}  \breve{x}_{i} \equiv \breve{l}_{1} \breve{x}_{1} + \breve{l}_{2} \breve{x}_{2} +\breve{l}_{3} \breve{x}_{3} = 0;
\end{equation}
the parameter $q$ is called, in analogy with ordinary finite projective planes, the order of  PR$_{\otimes}(2, q)$.
Let us find the total number of points/lines of PR$_{\otimes}(2, q)$. To this end, one first notes that from an algebraic point of view there
are two distinct kinds of them: I) those represented by the triples with at least one $\breve{x}_{i}/\breve{l}_{i}$ being a unit and II) those whose
representing $\breve{x}_{i}/\breve{l}_{i}$
are all zero-divisors, not all from the same ideal. It is then quite an easy exercise to see that PR$_{\otimes}(2, q)$ features
\begin{eqnarray}
\#^{{\rm (I)}}_{{\rm trip}} = \frac{\#_{{\rm t}}^{3} - \#_{{\rm z}}^{3}}{\#_{{\rm u}}}
&=& \frac{(q^{2})^{3} - (2q - 1)^{3}}{(q - 1)^{2}} =  \nonumber \\ &=& q^{4} + q^{2}(2q - 1) + (2q - 1)^2= (q^{2} + q + 1)^{2} - 6q
\end{eqnarray}
points/lines of the former type and
\begin{eqnarray}
\#^{{\rm (II)}}_{{\rm trip}} = \frac{\#_{{\rm z}}^{3} - \#_{{\rm s}}}{\#_{{\rm u}}} &=& \frac{(2q - 1)^{3} - (2q^{3} - 1)}{(q - 1)^{2}} = 6q
\end{eqnarray}
of the latter one; here $\#_{{\rm s}}$ stands for the number of distinct triples with all the entries in the same ideal.
Hence, its total point/line cardinality amounts to
\begin{equation}
\#_{{\rm trip}} = \#^{{\rm (I)}}_{{\rm trip}} + \#^{{\rm (II)}}_{{\rm trip}} = (q^{2} + q + 1)^{2}.
\end{equation}
Following the same chain of arguments/reasoning, but restricting to the classes of ordered couples instead, one finds that a line of
PR$_{\otimes}(2, q)$ is endowed with
\begin{equation}
\#_{{\rm coup}} = (q + 1)^{2}
\end{equation}
points and, dually, a point of PR$_{\otimes}(2, q)$ is the meet of the same number of lines.

Perhaps the most remarkable and fascinating feature of projective ring geometries is the fact that two distinct points/lines need not have
a unique connecting line/meeting point. More specifically, two distinct points/lines of PR$_{\otimes}(2, q)$ are called {\it neighbour} if they are
joined by/meet in at least {\it two} different lines/points; otherwise, they are called {\it distant} (or, by some authors, also
remote).\footnote{It is crucial to emphasize here that our concepts `neighbour' and/or `distant' are of a pure algebraic geometrical origin and
have {\it no}thing to do with the concept of metric --- the concept that does {\it not} exist {\it a priori} in any projective plane (space).}
Let us have a closer look at these interesting concepts. We shall pick up two different points of the plane, $A$ and $N$, where the former is
regarded as fixed and the latter as variable, and choose the coordinate
system in such a way that $A$ will be represented, without loss of generality, by the class
\begin{equation}
A:~  (1, 0, 0)
\end{equation}
whilst the representation of the other will be a generic one, i.e.
\begin{equation}
N:~  (\varrho a,  \varrho b, \varrho c).
\end{equation}
Our first task is to find out which constraints are to be imposed on $a$, $b$, and $c$ for $A$ and $N$ to be neighbours.
Clearly, any line ${\cal L}_{AN}$ passing through both $A$ and $N$ is given by
\begin{equation}
{\cal L}_{AN}:~  (0, \varsigma \beta, \varsigma \gamma)
\end{equation}
where
\begin{equation}
\beta b + \gamma c = 0.
\end{equation}
Next, we shall demonstrate that for $A$ and $N$ to be joined by more than one line, both $b$ and $c$ must be zero-divisors. For suppose that one
of the quantities, say $b$, is a unit; then from the last equation we get
\begin{equation}
\beta = - \gamma cb^{-1}
\end{equation}
which implies that
\begin{equation}
{\cal L}_{AN}:~  (0, -\varsigma \gamma cb^{-1}, \varsigma \gamma).
\end{equation}
Now,  $\gamma$ cannot be a zero-divisor because then all the entries in (14) would be zero-divisors of the same ideal; hence, it must be a unit which
means that
\begin{equation}
{\cal L}_{AN}:~  (0, -cb^{-1}, 1)
\end{equation}
which indeed represents, for fixed $a$ and $b$, just a {\it
single} class.  Returning to Eq.\,(12) yields that $b$ and $c$
belong to the same ideal. All in all, the coordinates of any
neighbour point $N$ to the point $A$ must be of the following two
forms
\begin{equation}
(\varrho a, \varrho g_{2} e, \varrho g_{3}e)~~{\rm or}~~(\varrho a, \varrho h_{2}(e - 1), \varrho h_{3}(e -1))
\end{equation}
where $a \in R_{\otimes}(q)$  and $g_{2}, g_{3}, h_{2}, h_{3} \in$
GF($q$), with the understanding that $a$, if being a zero-divisor,
is not from the same ideal as the remaining two entries and that
$g_{2}$ and $g_{3}$ or $h_{2}$ and $h_{3}$ do not vanish
simultaneously (which ensures that $N \neq A$). At this place a
natural question emerges: What is the cardinality of the
neighbourhood of a point of PR$_{\otimes}(2, q)$, that is, how
many distinct points $N$ are there? We shall address this question
in two steps according as $a$ is a unit, or a zero-divisor. In the
former case, taking $\varrho = a^{-1}$ brings Exps.\,(16) into
\begin{equation}
(1, g_{2}'e, g_{3}'e)~~{\rm or}~~(1, h_{2}'(e - 1), h_{3}'(e - 1))
\end{equation}
from which we infer that there are altogether 2($q^{2} - 1$) points of this kind. In the latter case, we have
\begin{equation}
(\varrho g_{1}(e - 1), \varrho g_{2}e, \varrho g_{3}e)~~{\rm or}~~(\varrho h_{1}e, \varrho h_{2}(e - 1), \varrho h_{3}(e -1))
\end{equation}
with $g_{1}, h_{1} \in$ GF($q$)$\setminus\{0\}$. Here, if  $g_{2}$, or  $h_{2}$, is non-zero, we can take $\varrho = - g_{1}^{- 1} (e - 1)
+ g_{2}^{-1}e$, or $\varrho = h_{1}^{-1}e -  h_{2}^{-1}(e - 1)$,  and reduce the last expressions into the forms
\begin{equation}
(e - 1, e, g_{3}'e)~~{\rm or}~~(e, e - 1, h_{3}'(e -1))
\end{equation}
which is the set of 2$q$ distinct points. If, finally, $g_{2}$
($h_{2}$) is zero (and, so,  $g_{3}$ ($h_{3}$) necessarily
non-zero), then the options $\varrho = - g_{1}^{-1}(e - 1) +
g_{3}^{-1}e$ ($\varrho = h_{1}^{-1}e - h_{3}^{-1}(e - 1)$)
transform Exps.\,(18) into
\begin{equation}
(e - 1, 0, e)~~{\rm or}~~(e, 0, e -1)
\end{equation}
which represent two different points. Hence, a point of PR$_{\otimes}(2, q)$ has altogether $2(q^{2} - 1) + 2q + 2 = 2q(q + 1)$ neighbour points.

\begin{figure}[ht]
\centerline{\includegraphics[width=6.4truecm,clip=]{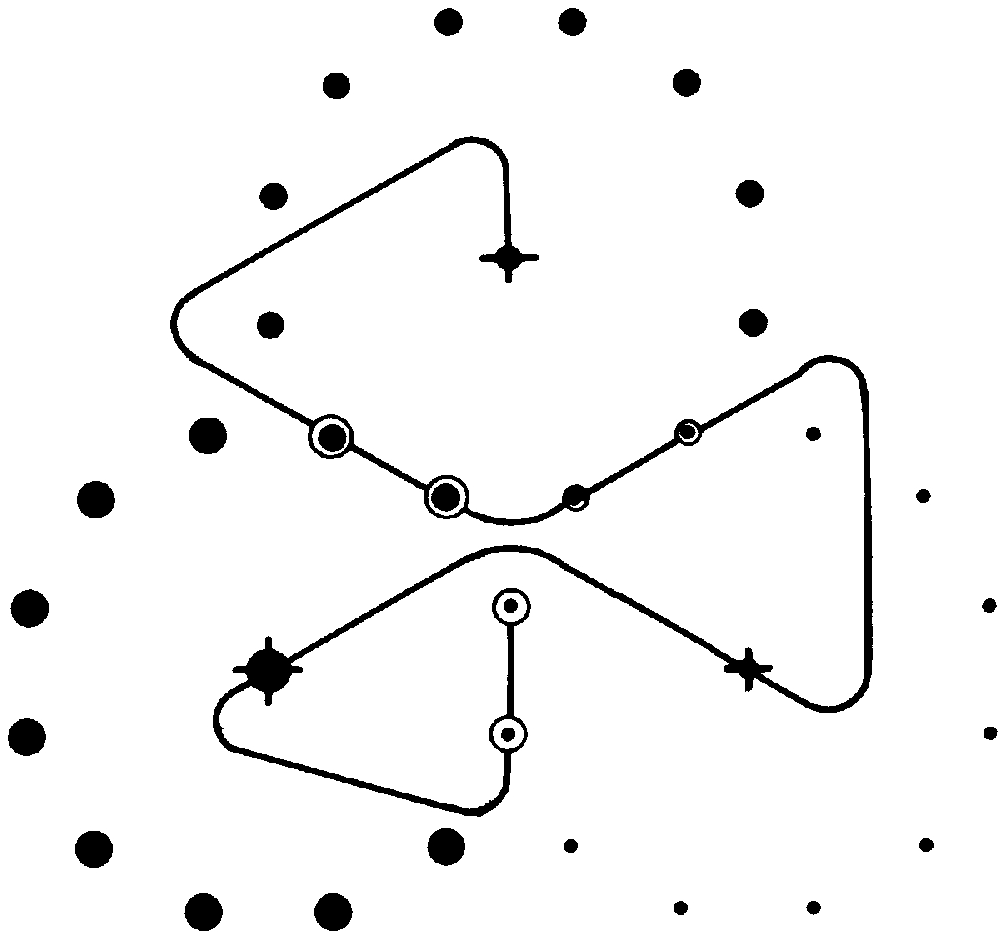}\includegraphics[width=5.7truecm,clip=]{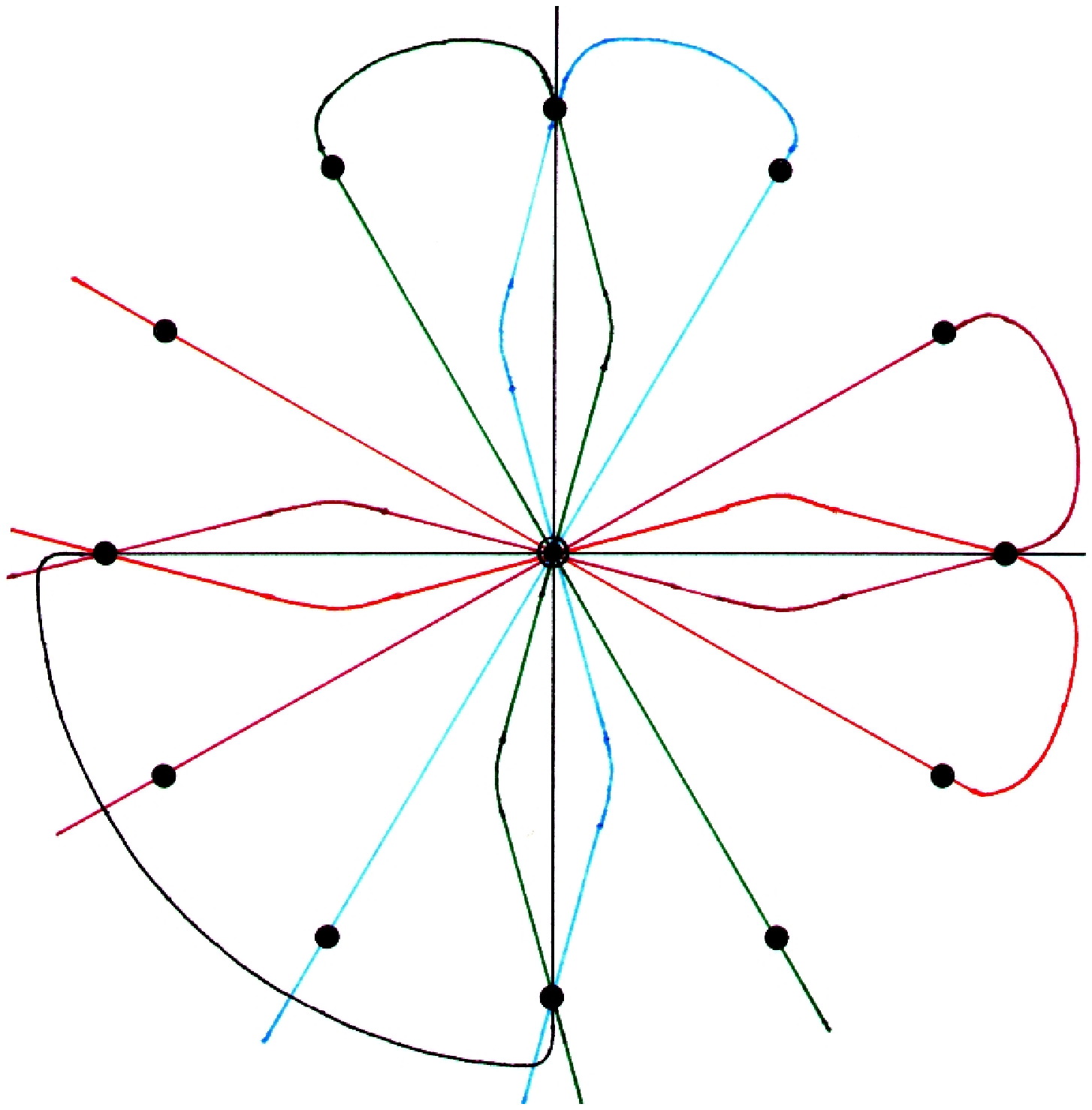}}
\caption{Some of the most salient properties of PR$_{\otimes}(2,
q$=2). {\it Left}: A line  (``bendy" curve) passes through the
common points of the neighbourhoods (sets of twelve points located
on different circles) of  any three mutually distant points
(crossed bullets, the centers of the circles) lying on it. {\it
Right}:  Each of lines (coloured curves) through a given point
(double circle) passes through four points of its neighbourhood.
This figure also exemplifies that two neighbour lines (e.g., the
blue and green ones) meet in three points and two neighbour points
(e.g., the double circle and the uppermost bullet) are joined by
three lines (blue, green and black in this case).}
\end{figure}

Next, let us consider a point $B$ that is {\it distant} to $A$; if we take, to facilitate our reasoning---yet not losing generality,
\begin{equation}
B:~  (0, 0, 1)
\end{equation}
then it is readily verified that ${\cal L}_{AB}:~  (0, 1, 0)$ is the only line connecting the two points. Employing here the above-introduced strategy,
the neighbourhood of the point $B$ is found to comprise $2(q^{2} - 1)$ points defined by
\begin{equation}
(g_{1}' e, g_{2}' e, 1)~~{\rm or}~~(h_{1}' (e - 1), h_{2}' (e - 1), 1),
\end{equation}
the 2$q$ ones of the form
\begin{equation}
(g_{1}' e, e, e - 1)~~{\rm or}~~(h_{1}' (e - 1), e - 1, e),
\end{equation}
and a couple represented by
\begin{equation}
(e, 0, e - 1)~~{\rm or}~~(e - 1, 0, e).
\end{equation}
Comparing these expressions with those of (17), (19) and (20), respectively, we find out that the two neighbourhoods, as the neigbourhoods
of {\it any} two distant points, are not disjoint but always (i.\,e., irrespective of the order of the plane) share {\it two} points (those defined by (20)/(24)
in our particular case). This means, algebraically, that the neighbour relation is not transitive and, so, it is not an equivalence relation
and, geometrically, that the neighbour classes (to a set of pairwise distant points) do not partition the plane. This important feature stems from
the fact that the ring $R_{\otimes}(q)$ is not local \cite{vk81},\,\cite{vk95}.
It can also be shown that the neighbourhoods of any three mutually distant points have never any point in common. Further, given a line and $q+1$
mutually distant points lying on it, the remaining points on the line are precisely those points in which the neighbourhoods of the chosen $q+1$ points
overlap; this
claim is easily substantiated by a direct cardinality check $q+1 + 2 \left( \begin{array}{c} q+1\\ 2 \end{array} \right) = q+1 + 2 (q+1)q/2 = (q+1)^{2}$.
Next, given a point, any line passing through the point is incident with $2q$ of its neighbour points. Finally, we mention that there are $q + 1$ lines
through  two distinct neighbour points and, dually, there are $q + 1$ points shared by two distinct neighbour lines. Fig.\,1 helps us visualise
some of these properties for the simplest case $q = 2$.

\section{Two Homomorphisms  PR$_{\otimes} (2, q) \mapsto $ PG$(2, q)$ and an\\ Algebraic Geometrical Complementarity Principle}
As already-mentioned, $R_{\otimes}(q)$ features two maximal ideals, Eqs.\,(2),\,(3), which implies the existence of two fundamental
homomorphisms,
\begin{equation}
\widehat{\pi}:~ R_{\otimes}(q) \mapsto \widehat{R}_{\otimes}(q) \equiv R_{\otimes}(q)/{\cal I}_{ \langle e \rangle} \cong {\rm GF}(q)
\end{equation}
and
\begin{equation}
\widetilde{\pi}:~ R_{\otimes}(q) \mapsto \widetilde{R}_{\otimes}(q) \equiv R_{\otimes}(q)/{\cal I}_{ \langle e-1 \rangle} \cong {\rm GF}(q),
\end{equation}
which induce two {\it complementary}, not-neighbour-preserving homomorphisms of  PR$_{\otimes} (2, q)$ into PG$(2, q)$, the ordinary (Desarguesian)
projective plane
of order $q$.\footnote{A general account of basic properties of a homomorphism between two projective ring planes/spaces can be found,
for example, in \cite{vk85} and \cite{vk95}, pp.\,1053--6.}
The nature of this complementarity is perhaps best seen on the behaviour of the neighbourhood of a point of  PR$_{\otimes} (2, q)$;
for applying
$\widehat{\pi}$ on Exps. (17), (19) and (20) yields, respectively,
\begin{equation}
(1, 0, 0)~~{\rm or}~~(1, -\widehat{h}_{2}', -\widehat{h}_{3}'),
\end{equation}
\begin{equation}
(1, 0, 0)~~{\rm or}~~(0, 1, \widehat{h}_{3}')
\end{equation}
and
\begin{equation}
(1, 0, 0)~~{\rm or}~~(0, 0, 1),
\end{equation}
whereas the action of $\widetilde{\pi}$ on the same expressions leads, respectively, to
\begin{equation}
(1, -\widetilde{g}_{2}', -\widetilde{g}_{3}')~~{\rm or}~~(1, 0, 0),
\end{equation}
\begin{equation}
(0, 1, \widetilde{g}_{3}')~~{\rm or}~~(1, 0, 0)
\end{equation}
and
\begin{equation}
(0, 0, 1)~~{\rm or}~~(1, 0, 0).
\end{equation}
From these mappings two important facts are readily discerned. First, under both the homomorphisms a {\it half} of the neighbour points merge with
the PG$(2, q)$ image of the point itself, while the other {\it half} of them go in a one-to-one correspondence to the remaining $(q^{2} + q + 1) - 1 = q(q+1)$
points of PG$(2, q)$;  for the simplest case, $q = 2$,  this feature is also illustrated in Fig.\,2. Second, the two sets play {\it reverse/complementary} roles
when switching from one homomorphism to the other; that is, the neighbours that merge together under one mapping  spread out
under the other mapping, and {\it vice versa}.

\begin{figure}[t]
\centerline{\includegraphics[width=11truecm,clip=]{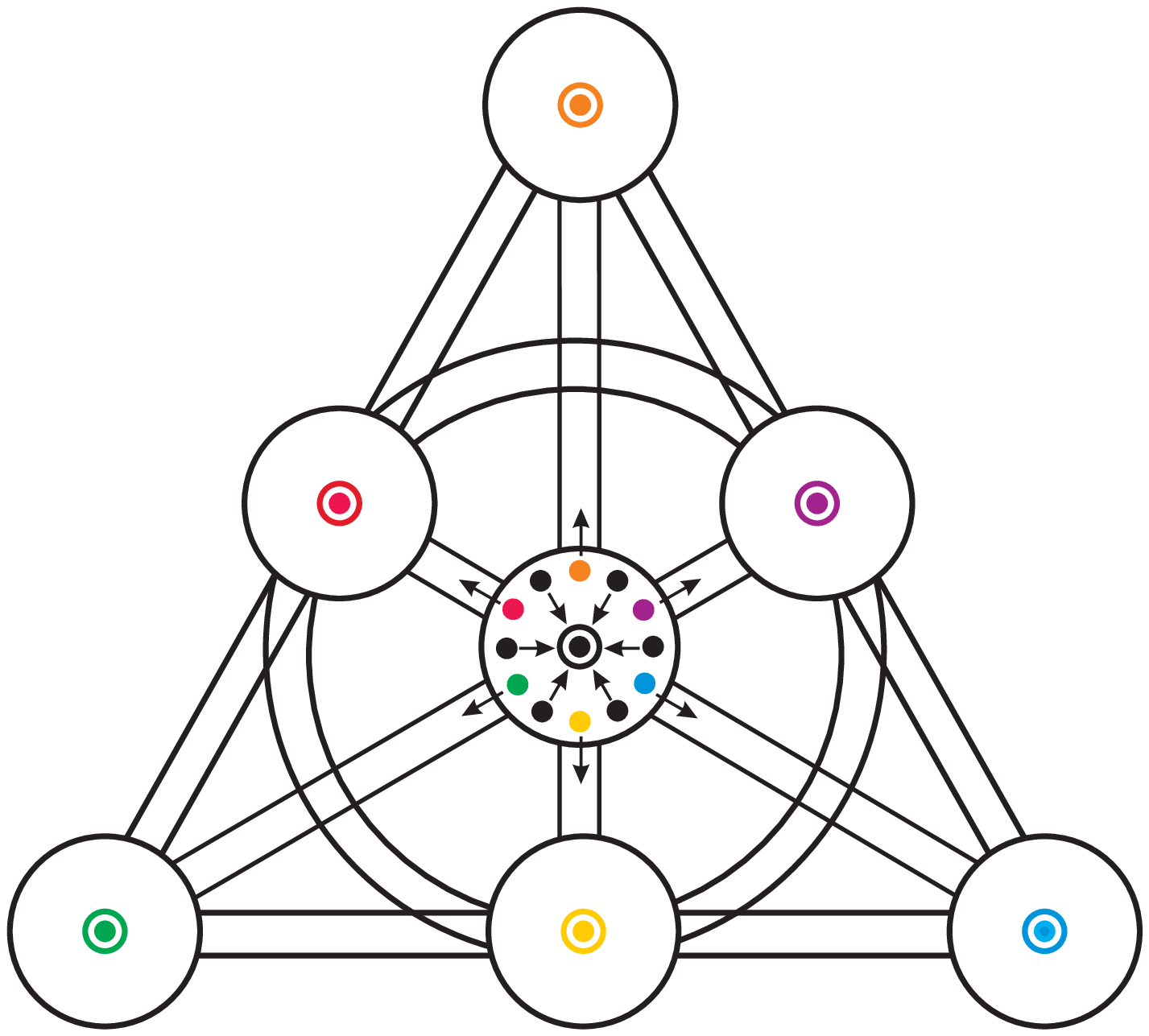}}
\caption{A schematic illustration of the properties of a homomorphism of PR$_{\otimes} (2, 2)$ into PG$(2, 2)$ on the behaviour of the
neighbourhood (twelve small points in the central circle) of a point of PR$_{\otimes} (2, 2)$ (the big black doubled circle). The points of the associated
(Fano) plane, PG$(2,2)$, are represented by seven big circles, six of its lines are drawn as pairs of line segments and the remaining line as a pair of
concentric circles. As indicated by small arrows, the six points  out of the twelve neighbours that are drawn black merge with the image of the
reference point itself, while those six points drawn in different colours are sent in a one-to-one, colour-matching fashion into the remaining
six points (the big coloured doubled circles) of PG$(2, 2)$.}
\end{figure}

Inherent to the structure of PR$_{\otimes} (2, q)$ is thus a remarkably simple, algebraic geometrical {\it principle of
complementarity}, which manifests itself in every geometrical object living in this plane. In order to get a deeper insight into its nature, let us
consider the points of a {\it line}. If we take, without loss of generality, the line to be the $(1, 0, 0)$ one, its points are the following ones: $q^{2}$
of them represented by the classes $(0, 1, r)$, where $r \in R_{\otimes}(q)$, $2q - 1$ defined by $(0, d, 1)$, $d$ being a zero-divisor, and the rest
of the form $(0, d_{I}, d_{II})$, with both $d_{I}$ and $d_{II}$ being zero-divisors not of the same ideal. More explicitly, the points of the first set
comprise the points $(0, 1, u)$, $u$ being a unit, the point $(0, 1, 0)$ and $2(q-1)$ points of the form
\begin{equation}
(0, 1, ge)~~{\rm or}~~(0, 1, h(e-1)),
\end{equation}
those of the second set feature the point $(0, 0, 1)$ and  $2(q-1)$ points defined by
\begin{equation}
(0, ge, 1)~~{\rm or}~~(0, h(e-1), 1),
\end{equation}
and the remaining two points are
\begin{equation}
(0, e, e-1)~~{\rm or}~~(0, e-1, e).
\end{equation}
If these expressions are subject to the two homomorphisms, the former gives
\begin{equation}
(0, 1, 0)~~{\rm or}~~(0, 1, -\widehat{h}),
\end{equation}
\begin{equation}
(0, 0, 1)~~{\rm or}~~(0, -\widehat{h}, 1),
\end{equation}
\begin{equation}
(0, 0, 1)~~{\rm or}~~(0, 1, 0),
\end{equation}
respectively,  whereas the latter produces
\begin{equation}
(0, 1, -\widetilde{g})~~{\rm or}~~(0, 1, 0),
\end{equation}
\begin{equation}
(0, -\widetilde{g}, 1)~~{\rm or}~~(0, 0, 1),
\end{equation}
\begin{equation}
(0, 1, 0)~~{\rm or}~~(0, 0, 1),
\end{equation}
respectively. Again, the complementarity of the two mappings is well pronounced.

The structure of PR$_{\otimes} (2, q)$, embodied in the properties of $R_{\otimes}(q)$, thus admits a remarkable complementary
interpretation/description in terms of the properties of the associated ordinary projective plane, PG$(2, q)$, and {\it either} of the mappings,
$\widehat{\pi}$ or $\widetilde{\pi}$. From the above-given examples it is obvious that either of the representations is partial and insufficient
by itself; separately neither of them fully grasps the properties of  PR$_{\otimes} (2, q)$, they do that only when taken together. Hence,
PR$_{\otimes} (2, q)$ can serve as an elementary, algebraic geometrical expression of the principle of complementarity! Let us highlight
its potential domains of application.

\section{Possible Applications of the Geometry/Configuration}
As it is very well known, the principle of complementarity was
first suggested by Niels Bohr \cite{bohr} in an attempt to
circumvent severe conceptual problems at the advent of quantum
mechanics.  Quantum theory is thus the first domain when we should
look for possible applications of the geometry of PR$_{\otimes}
(2, q)$. This view, in  fact, gets strong support from our recent
work \cite{spl} where, as already mentioned, we have shown that a
closely related class of finite ring planes, projective Hjelmslev
planes, provide important clues as per probing the structure of
finite-dimensional Hilbert spaces. Another promising physical
implementation of PR$_{\otimes} (2, q)$ turns out to be the
concept of  an abstract primordial prespace, stemming from the V\"
axj\" o interpretation of quantum mechanics \cite{sk}. Here, three
different fundamental ``sectors" of the prespace are hypothetised
to exist, namely classical, semiclassical and pure quantum
according as the coordinatizing ring of the corresponding
projective plane --- taken in the simplest case to be a quadratic
extension of a Galois field --- is a field, a local ring or a ring
with two maximal ideals (of zero-divisors), respectively. In this
framework, the phenomena like quantum non-locality and quantum
entanglement may well find their alternative explanations in terms
of neighbour/distant relations and ring-induced homomorphisms may
provide natural pairwise couplings between the sectors of the
prespace. It has as well been suspected that finite ring planes
and related combinatorial concepts/configurations (e.g.,
zero-divisor graphs) may help overcome some technical problems
when describing certain highly complex macro-molecular systems
\cite{sp}, as their subsystems often exhibit traits of puzzling
dual, complementary behaviour. In this context, it is worth
mentioning some intriguing parallels with El Naschie's 
``Cantorian" fractal approach to quantum mechanics (see, e.g.,
\cite{eln04}--\cite{eln06}), for this remarkable theory --- whose foundations rest on several well-established branches of mathematics including
topology, set theory, algebraic geometry and group/number theory --- has also been
recognized to be, in principle, extendible over algebras/number
systems having {\it zero-divisors} \cite{cz}.

If we go beyond physics, we are likely to find even more intriguing phenomena where  PR$_{\otimes} (2, q)$ could be employed.
One is encountered in the field
of neurophysiology and it is called perceptual rivalry (see, e.\,g., \cite{pet}). It is, roughly speaking, the situation when there are two
different, competing with each other interpretations that our brain can make of some sensory event. It is a sort of oscillation,
due to interhemispheric switching, of conscious experience between two complementary modes despite univarying sensory input.
Its particular case is the so-called binocular rivalry, i.\,e. alternating  perceptual states that occur when the images seen by both the eyes
are too different to be fused into a single percept. Here, we surmise the two homomorphims of PR$_{\otimes} (2, q)$ into PG$(2, q)$ to be
capable of qualitatively
underlying the activities in the two hemispheres of the brain, with the associated Galois plane standing for a mediator (``switch") between them.

Clearly, a lot of --- mostly conceptual --- work is to be done along many lines of inquiry in order to see whether the exciting prospects implicit in these
conjectures are real or merely illusory. The structure of PR$_{\otimes} (2, q)$ is, however, so enchanting that such work is certainly worth pursuing and
our recent closely-related investigations \cite{psk} indeed seem to justify this feeling.\\ \\ \\
\noindent
\Large
{\bf Acknowledgements}
\normalsize

\vspace*{.2cm}
\noindent
The first author thanks Mr. Pavol Bend\' {\i}k and Dr. Petr Pracna for a careful drawing of the figures and Dr. Richard Kom\v z\' {\i}k 
for a computer-related
assistance. This work was partially supported by the
Science and Technology Assistance Agency under the contract $\#$
APVT--51--012704, the VEGA project $\#$ 2/6070/26 (both from
Slovak Republic) and by the trans-national ECO-NET project $\#$
12651NJ ``Geometries Over Finite Rings and the Properties of
Mutually Unbiased Bases" (France).

\end{document}